\documentstyle{brkstyle}

\newcommand{\decale}{\hspace*{3ex}}
\newcommand{\ld}{,\ldots,}

\newcommand{\pl}{\partial}
\newcommand{\ol}{\overline}
\newcommand{\olp}{\ol\pl}
\newcommand{\old}{\ol D}

\newcommand{\sub}{\subset\!\subset}
\newcommand{\ww}{\wedge\cdots\wedge}

\newcommand{\re}{\mathop{\rm Re}\nolimits}

\renewcommand{\ln}{\ell\!n}
\newcommand{\sgn}{\mathop{\rm sgn}\nolimits}

\begin{document}

\begin{center}
\title{\bf  Optimal Regularity for $\overline\partial\sb {\rm b}$ on $CR$ manifolds}\\
{\em {\small by Moulay Youssef Barkatou
}}
\end{center}

\Names{   Barkatou}
\Shorttitle{  Optimal Regularity...}  
\Subjclass{32F20-32F10-32F40.}
\Keywords{CR manifold, tangential Cauchy-Riemann equations, $q$-convexity}

\maketitle

\begin{abstract}
In this paper a new  explicit integral formula is derived for
solutions of the tangential Cauchy-Riemann equations on $CR$ $q$-
concave manifolds and optimal estimates in the Lipschitz norms are obtained.  

\end{abstract}

\setcounter{section}{-1}

\newsection{\hspace*{-1.5ex}Introduction} 
 The aim of this paper is to prove the following theorem: 
\paragraph{{\sc  Theorem 0.1}}{\em Let $M$ be a $q$-concave $CR$
 generic submanifold (cf.sect 1.2) of codimension $k$ and of class
 ${\cal C}^{2+\ell}$(resp. ${\cal C}^{3+\ell}$) in $\C^n$ ($\ell\ge
 0$) and $ z_0$ a point  in $M$. Then
there exist an open neighborhood $M_0\subseteq M$ of $z_0$ and kernels ${\cal R}_r(\zeta,z)$, for $r=0\ld q-1,n-k-q\ld n-k$, with the following properties,

\decale{(i)} For every domain $\Omega\sub M_0$ with piecewise ${\cal
  C}^1$ boundary and every   ${\cal C}^1$ $(0,r)$-form $f$ on
$\overline \Omega$ ( $0\le r\le q-1$ or $n-k-q+1\le r\le n-k$ ), we have
$$f=\overline\partial\sb {\rm b}\int_\Omega f\wedge{\cal R}_{r-1} - \int_\Omega \overline\partial\sb {\rm b}f\wedge{\cal R}_r +\int_{{\rm b}\Omega} f\wedge{\cal R}_r$$
on $\Omega$.

\decale{(ii)} For every open  set $\Omega\sub M_0$ the integral
operator $\int_\Omega \cdot\wedge{\cal R}_r$ is a bounded linear
operator from  ${\cal C}^{\ell}_{0,r+1}(\Omega)\cap L^{\infty}(\Omega)$ to ${\cal C}^{\ell+{1\over 2}}_{0,r}(\Omega)$ for $r\ge n-k-q$ (resp. $r\le q-1)$}.
 
\vspace*{0.5cm}

\noindent Theorem 0.1 has the following corollary
 
\paragraph{{\sc  Corollary 0.2}}{\em Let $M$ be a $1$-concave $CR$ generic ${\cal C}^{3+\ell}$-submanifold of a complex manifold. Let $T$ be a distribution of order $\ell$ on $M$. If $\overline\partial\sb {\rm b}T$ is defined by a ${\cal C}^\ell$ $(0,1)$-form on $M$ then $T$ is defined by a ${\cal C}^{\ell+{1\over 2}}$- function on $M$.}

\vspace*{0.5cm}

 For a proof of corollary 0.2, we refer to the proof of
theorem 4.1.6 in  \cite{[Ba1]}. The interest of this Corollary lies in the fact that under the hypothesis of 1-concavity the tangential Cauchy-Riemann equation for (0,1)-currents  cannot be solved locally (see \cite{[A/F/N]}).

Theorem 0.1 and Corollary 0.2 essentially improve the results of
Airapetjan and Henkin \cite{[He1]}, \cite {[Ai/He1]}, \cite {[Ai/He2]}
and also of  the author in \cite{[Ba1]} where homotopy formulas were
obtained with less explicit kernels giving almost optimal but not optimal estimates.

The study of the tangential Cauchy-Riemann equations with the use of
explicit integral formulas was initiated by Henkin \cite{[He1]}. For
further references and results on $CR$ manifolds we refer the reader
to the survey of Henkin \cite{[He2]}, the memoir of Tr\`eves \cite{[Tr]} and the book of Boggess \cite{[Bo]}.

It is known that a fundamental solution for the $\overline\partial\sb
{\rm b}$ operator on certain hypersurfaces (see Henkin \cite{[He3]},
Harvey-Polking \cite{[H/P]}, Boggess-Shaw \cite{[Bo/Sh]}, Fischer-Leiterer \cite{[Fi/Le]}) can be constructed as the jump of two kernels, obtained by applying to the usual Bochner-Martinelli-Koppelman kernel (BMK kernel) in $\C^n$, a solution operator for $\overline\partial$, once on the left and once on the right hand side of the hypersurface.

Solutions  for such equations  can be given by applying the
generalized Koppelman (cf.section 1.3) to the BMK section and the
barrier functions (cf.section 1.4) of the hypersurface as was done  in
\cite{[He3]}, \cite{[H/P]}, \cite{[Bo/Sh]} and \cite{[Bo]} or by using a homotopy
operator for $\overline\partial$ of Grauert-Lieb-Henkin type as was
achieved in \cite{[Fi/Le]}.

Inspired by the definition of a hyperfunction of several variables,
the present author generalized in \cite{[Ba1]} the construction of
Fischer-Leiterer \cite{[Fi/Le]} to higher codimensional $CR$
submanifolds by solving with estimates up to the boundary some $\overline\partial$ equations on certain wedges attached to such manifolds with the use of $\overline\partial$ homotopy operators from  \cite{[La/Le1]} and \cite{[La/Le2]}.

In this paper we shall show that such equations can also be solved up
to some error terms by using the Koppelman lemma (see (2.2)) and the
key idea in this work is to "deform" via this lemma those terms into
ones with vanishing coefficients for some bidegrees (see lemmas 2.2
and 2.3), the strict $q$-convexity plays here an important role.

 We shall give two fundamental solutions to the tangential
 Cauchy-Riemann complex. The first one (cf.sect 2.1) does not yield
 sharp estimates for the solutions of $\overline\partial\sb {\rm b}$
 (when $k>1$) but is  a "necessary" step to construct the second one
 (cf. sect 2.2)  corresponding to kernels ${\cal R}_r$. To derive the
 latter fundamental solution from the former, we shall use an idea of
 Henkin \cite{[He3]}.

\vspace*{0.3cm}

In \cite {[Fi]}~B.Fischer proved Theorem 0.1 and Corollary 0.2 for
hypersurfaces by using a version of the first fundamental solution
which was suggested to him by I.Lieb and J.Michel.

Recently, Polyakov \cite {[Po]} proved  sharp estimates for global solutions of $\overline\partial\sb {\rm b}$ on $q$-concave $CR$ manifolds, in Lipschitz spaces of Stein \cite {[St]}.  
\paragraph{{\sc \bf  Polyakov's theorem.}}{\em
Let $M$ be a $q$-concave $CR$ generic ${\cal C}^4$-submanifold in $\C^n$ with $q\ge 2$ and let $M'$ be a relatively compact open subset of $M$. Then for any $r=1\ld q-1$ there exist linear operators
$$R_r:L^s_{(0,r)}(M)\rightarrow\Gamma^{s,1}_{(0,r-1)}(M) {~and~}  H_r:L^s_{(0,r)}(M)\rightarrow L^s_{(0,r)}(M)$$
such that for any $s\in[1,\infty]$ $R_r$ is bounded and $H_r$ is compact and such that for any differential form $f\in {\cal C}^{\infty}_{(0,r)}(M)$ the following equality:
$$f(z)=\overline\partial\sb {\rm b}R_r(f)(z)+R_{r+1}(\olp\sb {\rm b}f)(z)+H_r(f)(z)$$
holds for $z\in M'.$}

\vspace*{0.5cm}

Our method is quite different from that of Polyakov, and it is not clear how one can get an analogous result to Corollary 0.2 from Polyakov's theorem.

\vspace*{0.3cm}

This paper is organized as follows. In section 1.2 we give the definition of a $q$-concave $CR$  manifold and we define the $\overline\partial\sb {\rm b}$ operator. In section 1.3 we  recall the generalized Koppelman lemma which plays a key role in the construction of our kernels. In section 1.4 we  recall the construction of a barrier function  and  a Leray map for a hypersurface at a point where the Levi form has some positive eigenvalues. In section 1.5 we state some elementary facts from Algebraic Topology, which we shall use later. Section 2.1 is devoted to the construction of our first fundamental solution. In section 2.2 we construct our second fundamental solution and in section 3 we prove estimates for our kernels.

\section{\hspace*{-1.5ex} Preliminaries and notations}

\paragraph{1.1.}Let $X$ be a complex manifold, and M a real
submanifold of X .

Let $f$ be a differential form of degree m defined on a domain  $D\subseteq
M$. Then we denote by $\|f(z)\|, z\in D$, the Riemannian norm of $f$ at $z$
(cf.\cite {[He/Le2]}, section 0.4), and we set

  $$ \|f\|_{0,D} = \sup_{z\in D} \|f(z)\|$$

\noindent and

$$\|f\|_{\alpha ,D}=\|f\|_{0,D} + \sup_{{z,\zeta }\atop{z\ne\zeta
}} {\|f(z)-f(\zeta )\|\over |\zeta -z|^\alpha }$$
for $0<\alpha <1$ .

If $0<\alpha <1$ , then a form on $D$ is called $\alpha$-H\"older continuous on $D$ if
                $$ \|f\|_{\alpha ,K}\le \infty .$$ for all compact sets $K\subseteq D$.

If $\ell$ is a non-negative integer and $0<\alpha <1$ , then we say  f is a ${\cal C}^{\ell+\alpha}$ form on $D$ if $f$ is of class ${\cal C}^{\ell}$ and all
 derivatives of order $\le \ell$ of $f$ are $\alpha$-H\"older
 continuous on $D$. $L^{\infty}(D)$ denotes the space of all bounded
 forms on $D$.

Throughout this paper $C$ will denote a positive constant which is
independent of the variables and the functions. The constant $C$ used in different places
may have different values there.

\paragraph{1.2.}Let $M$ be a real submanifold of class ${\cal C}^2$ in  $\C^n$ defined by
$$M=\{z\in\Omega; \rho _1(z)=\cdots = \rho _k(z)=0\}\quad 1\le k\le n\eqno
(1.1)$$
where $\Omega $ is an open subset of $\C^n$ and the functions    $\rho _\nu $, $1\le\nu \le k$, are real-valued functions of class ${\cal C}^2$ on  $\Omega $   with the property
$d\rho _1(z)\ww d\rho _k(z)\ne 0$ for each $z\in M$. 

We  denote by  $T^\C_z(M)$ the complex tangent space to $M$ at the point $z\in M$ i.e.,   

$$T^\C_z(M)=\{\zeta \in\C^n / \sum^n_{j=1} {\pl\rho _\nu \over \pl z_j} (z)
\zeta _j=0, \nu =1\ld k\}.$$
We have $\dim_\C T^\C_z(M)\ge n-k$. The submanifold  $M$ is called a Cauchy-Riemann manifold ($CR$-manifold) if the number $\dim_\C T^\C_z(M)$ does not depend on the point $z\in M$. $M$ is said to be  $CR$
generic  if for every $z\in M$, $\dim_\C T^\C_z(M)=n-k$  , this is equivalent to~:
$$\olp\rho _1\wedge\olp\rho _2\ww\olp\rho _k\ne
0 \hbox{~on~} M.\eqno (1.2)$$

If $M$ is $CR$ generic , we call  $M$  $q$-concave,
$\displaystyle 0\le q \le {n-k\over 2}$, if for each $z\in M$ and every
$x\in\R^k\setminus\{0\}$ the following hermitian form on $T^\C_z(M)$
$$\sum_{\alpha ,\beta } {\pl^2\rho _x\over \pl z_\alpha \pl\ol z_\beta }(z)
\zeta _\alpha \ol\zeta _\beta ,\hbox{~where~} \rho _x=x_1\rho _1+\cdots x_k \rho
_k\ $$
has at least  $q$ negative eigenvalues.

If $M$ is $CR$ generic then we denote by ${\cal C}^s_{p,r}(M)$ the space of differential forms of type $(p,r)$ on $M$ which are of class ${\cal C}^s$. Here, two forms $f$ and $g$ in ${\cal C}^s_{p,r}(M)$ are considered to be equal if and only if for each form $\varphi\in  {\cal C}^\infty_{n-p,n-k-r}(\Omega)$ of compact support, we have   
$$ \int_M f\wedge \varphi = \int_M g\wedge \varphi .$$
We denote by ${\cal L}^{(-s)}_{p,r}(M)$ the dual space to ${\cal C}^s_{n-p,n-k-r}(M)$.

We define the tangential Cauchy-Riemann operator on forms in ${\cal L}^{(-s)}_{0,r}(M)$ as follows. If $u\in {\cal C}^s_{0,r}(M)$, $s\ge 1$, then $u$ can be extended to a smooth form $ \tilde u\in{\cal C}^s_{0,r}(\Omega)$ and we may set
$$\overline\partial\sb {\rm b}u :=\overline\partial \tilde u|_M$$     
 It follows from the condition for equality of forms on $M$ that this definition does not depend on the choice of the extended form $\tilde u$. In general, for for forms $u\in{\cal L}^{(-s)}_{0,r-1}(M)$ and $f\in{\cal L}^{(-s)}_{0,r}(M)$, by definition,
$$\overline\partial\sb {\rm b}u =f$$
will mean that for each form $\varphi\in  {\cal C}^\infty_{n-p,n-k-r}(\Omega)$ of compact support we have
 $$ \int_M f\wedge \varphi ={(-1)^ r} \int_M u\wedge \overline\partial
 \varphi .$$

\paragraph{ 1.3. The generalized Koppelman lemma.}

In this section we recall a formal identity ( the generalized
Koppelman lemma) which is essential for the construction of our
kernels. The exterior calculus we use here was developed by Harvey and
Polking in \cite{[H/P]}.

Let $V$ be an open set of $\C^n \times \C^n$. Suppose $G:V\rightarrow\C^n$ is a ${\cal C} ^1$ map. We write

 $$ G(\zeta,z)=(g_1(\zeta,z),\ldots,g_n(\zeta,z))$$
and we use the following notation
$$G(\zeta,z).(\zeta-z) = \sum^n_{j=1}g_j(\zeta,z)(\zeta_j-z_j)$$   
$$G(\zeta,z).d(\zeta-z) = \sum^n_{j=1}g_j(\zeta,z)d(\zeta_j-z_j)$$
$${\overline\partial_{\zeta,z} G}(\zeta,z).d(\zeta-z) = \sum^n_{j=1}{\overline\partial_{\zeta,z} g_j}(\zeta,z)d(\zeta_j-z_j )$$
where $\overline\partial_{\zeta,z}={\overline\partial_\zeta}
+{\overline\partial_z}$.

\vspace*{0.3cm}

We define the Cauchy-Fantappie form $\omega^G$ by
$$\omega^G = {G(\zeta,z).d(\zeta-z)\over G(\zeta,z).(\zeta-z)}$$
on the set where $G(\zeta,z).(\zeta-z)\ne 0$.
  
\vspace*{0.3cm}

Given m such maps, $G^j$, $1\leq j\leq m$, we define the kernel
$$ {\Omega (G^1,\ldots,G^m)}=\omega^{G^1}\wedge\ldots\wedge\omega^{G^m}\wedge \sum_{\alpha_1+\ldots+\alpha_m=n-m}{(\overline\partial_{\zeta,z}\omega^{G^1})}^{\alpha_1}\wedge\ldots\wedge{(\overline\partial_{\zeta,z}\omega^{G^m})}^{\alpha_m}$$
on the set where all the denominators are nonzero.
\paragraph{{\sc  Lemme 1.1.}}{\em{\rm  (The generalized Koppelmann lemma)} $$\overline\partial_{\zeta,z}{\Omega(G^1,\ldots,G^m)}=\sum^m_{j=1}{(-1)^j}{\Omega(G^1,\ldots,\hat{G^j},\ldots,G^m)}$$
on the set where the denominators are nonzero, the symbol $\hat{G^j}$ means that the term $G^j$ is deleted.}

\vspace*{0.5cm}
    
The following lemma  is useful for the estimation of the kernel defined above.
\paragraph{{\sc  Lemme 1.2.}}{\em

For $k\ge 0$  
$$\omega^G \wedge (\overline\partial_{\zeta,z}\omega^G)^k
={G(\zeta,z).d(\zeta-z)\over G(\zeta,z).(\zeta-z))}\wedge
\Big({{\overline\partial_{\zeta,z}G.d(\zeta-z)\over
 G(\zeta,z).(\zeta-z)}} \Big)^k.$$}

\vspace*{0.5cm}

For a proof of these two lemmas we refer the reader to \cite{[H/P]} or \cite {[Bo]}.
\paragraph{{\sc  Remark.}}{\em
When  $G=\overline {\zeta - z}$ , we see from Lemma 1.2 that $\Omega (G)$
 is the classical Martinelli-Bochner Koppelman kernel in $\C^n$.} 

\paragraph{ 1.4. Barrier function.}

In this section, we shall construct a barrier  function for a hypersurface at a point where the Levi form has some positive eigenvalues.

For a detailed proof of what follows we refer the reader to sect.3 in \cite {[La/Le1]}.

Let $H$ be an oriented real hypersurface of class ${\cal C}^2$ in  $\C^n$ defined by
$$H=\{z\in\Omega; \rho (z)=0\} $$
where $\Omega $ is an open subset of $\C^n$ and $\rho  $ is a real-valued function of class ${\cal C}^2$ on  $\Omega $ with $d\rho (z)\ne 0$ for each $z\in H$.

Denote by  $F(\cdot,\zeta )$ the Levi polynomial of   $\rho $ at a point $\zeta \in\Omega $,
i.e.

$$F (z,\zeta )=2\sum^n_{j=1} {\pl\rho (\zeta )\over \pl\zeta _j} (\zeta
_j - z_j) - \sum^n_{j,k=1} {\pl^2\rho (\zeta )\over \pl\zeta _j \pl\zeta _k}
(\zeta _j - z_j)(\zeta _k - z_k)$$
$\zeta \in\Omega , z\in\C^n$.
 
Let $z^0\in H$ and  
 $T$ be the largest vector subspace of $\C^n$ such that the Levi form of $\rho$ at $z\sp0$ is positive definite on $T$. Set ${\bf dim T= d}$ and suppose $d\ge 1$.

Denote by $P$ the orthogonal projection from $C^n$ onto $T$, and set $Q=I-P$.
Then it follows from Taylor's theorem that there exist a number $R$ and two  positives constants $A$ and $\alpha$ such that the following holds:
$$\re F(z,\zeta )\ge \rho (\zeta ) - \rho (z) + \alpha |\zeta -z|^2 - A |Q (\zeta -z)|^2\eqno (1.3)$$
 for $|z^0-\zeta|\le R$ and $|z^0-z|\le R$.

 Since $\rho$ is of class ${\cal C}^2$ on $\Omega$ ,  We can find ${\cal C}^\infty$ functions $a^{kj}( k,j=1\ld n)$ on a  neighborhood $U$ of $z^0$  such that
$$\left\vert
a^{kj}(\zeta )-{\pl^2 \rho  (\zeta )\over \pl\zeta _k\pl\zeta
_j}\right\vert < {\alpha \over 2n^2}$$
for all $\zeta \in U$. And then we have 
$$\left\vert
\sum^n_{k,j=1} \Big(a^{kj}(\zeta )-{\pl^2\rho  (\zeta )\over \pl\zeta
_k \pl\zeta _j}\Big) t_k t_j\right\vert\le {\alpha \over 2} |t|^2$$
for all $\zeta \in U$ and $t\in\C^n$. Set
$$\widetilde{F}(\zeta,z )=2 \sum^n_{j=1} {\pl \rho 
(\zeta )\over \pl\zeta j} (\zeta_j - z_j) - \sum^n_{k,j=1} a^{kj}
(\zeta)(\zeta _k-z_k) (\zeta _j-z_j)$$
for $(z,\zeta)\in\C^n \times U$.Then it follows from $(1.3)$ that 
$$\re\widetilde{F}(\zeta,z )\ge \rho (\zeta ) - \rho (z) + {\alpha \over 2} |\zeta -z|^2 - A |Q(\zeta
-z)|^2\eqno (1.4) $$
for $|z^0-\zeta|\le R$ and $|z^0-z|\le R$.

 Denote by $Q_{kj}$ the entries of the matrix
$Q$ i.e
$$Q=\big(Q_{kj}\big)^n_{k,j=1}\quad (k=\hbox{~ column index}).$$
We set for $(z,\zeta)\in\C^n \times U$
$$\left\{\displaystyle
\begin{array}{rl}
  g_j (\zeta,z) & = 2{\pl \rho (\zeta
)\over \pl\zeta j} - \sum^n_{k=1} a^{kj} (\zeta ) (\zeta _k - z_k) + A
\sum^n_{k=1} \ol{Q_{kj}(\zeta _k - z_k)}\\
G(\zeta,z) & = \big(g_1(\zeta,z),\ld, g_n(\zeta,z) \big)\\
\Phi (\zeta,z) & = G(\zeta,z).(\zeta-z).
\end{array}\right.$$
Since $Q$ is an orthogonal projection, then we have
$$\Phi (\zeta,z)= \widetilde{F}(\zeta,z)+ A |Q(\zeta-z)|^2$$
hence it follows from $(1.4)$ that
$$\re\Phi(\zeta,z) \ge  \rho (\zeta ) - \rho (z) + {\alpha
\over 2} |\zeta -z|^2\eqno (1.5)$$
for $|z^0-\zeta|\le R$ and $|z^0-z|\le R$.

\vspace*{0.3cm}

$G$ is called a {\bf Leray map} and
$\Phi$ is called a {\bf barrier function} of $H$(or $\rho$)  at $z^0$.
\paragraph{\sc Definition.}{\em A map $f$ defined on some complex manifold $X$ will be called $k$-holomorphic if, for each point $\xi\in X$, there exist holomorphic coordinates $h_1 \ld h_k$ in a neighborhood of $\xi$ such that $f$ is holomorphic with respect to $h_1\ld h_k$.}
\paragraph{\sc Lemma 1.3.}{\em 
For every fixed $ \zeta \in U$, the map $G(\zeta,z)$ and the function $\Phi$, defined above, 
are $d$-holomorphic in $z\in \C^n$.}

\vspace*{0.5cm}

{\sl Proof}. Choose complex linear coordinates $h_1 \ld h_n$ on $\C^n$ with
$$\{z\in \C^n :Q(z)=0\}=\{z\in \C^n :h_{d+1}(z)=\cdots=h_n(z)=0\}.$$
Then the map $\C^n\ni z\rightarrow\overline{Q(\zeta-z)}$ is independent of
$h_1\ld h_d$. This implies that $G(\zeta,.)$ is complex linear with respect to
 $h_1\ld h_d$, and $\Phi(\zeta,.)$ is quadratic complex polynomial with respect to $h_1\ld h_d$.\hfill$\Box$
\paragraph{{1.5. Some Algebraic Topology.}}Let $N$ be a positive integer.Then   we call  p-simplex,  $1\le p\le N$, every collection of p  lineary independent vectors in $\R^N$.

We define $S_p$ as the set of all  finite formal  linear combinations, with integer coefficients, of $p$-simplices.

Let $\sigma =[a_1\ld a_p]$ be a $p$-simplex,  then we set
$$\pl_j \sigma =[a_1 \ld \hat{a_j} \ld a_p]$$ for $1\le j\le p$ and
$$\pl\sigma =\sum^p_{j=1} {(-1)}^j \pl_j \sigma$$
(this definition holds also for any collection of $p$ vectors). 
If $1\le j_1\le p \ldots 1\le j_r\le p-r$,  we define
$$\pl^r_{j_r \cdots j_1}\sigma = \pl_{j_r}(\pl^{r-1}_{j_{r-1} \cdots j_1}\sigma)$$
where $\pl^1_j\sigma =\pl_j\sigma $.

All of these operations can be extended by linearity to $S_p$.

 If $\sigma$ is a $p$-simplex defined as above then we define the barycenter of $\sigma$ by
$$b(\sigma) = {1\over p} \sum^p_{j=1}a_j.$$ 
Now we define the first barycentric subdivision of $\sigma$
by the following
$$ sd(\sigma) =(-1)^{p+1}\sum_{j_1 \ld j_{p-1}\atop {1\le j_i\le p-i+1}}   { (-1)^{j_1+\cdots +j_{p-1}}\Big[b(\sigma),b(\pl_{j_1}\sigma) \ld b(\pl ^{p-1}_{j_{p-1}\cdots j_1} \sigma)\Big ]}.$$
By linearity we can  also define the first barycentric Subdivision of any element of $S_p$.

It is easy to see that

\paragraph{{\sc Lemma 1.4.}}{\em
If $\sigma$ is an element of $S_p$, then
$$ sd(\pl\sigma) = \pl sd(\sigma). $$}
\vspace*{0.3cm}

The barycentric subdivision of higher order of an element $\sigma$ of $S_p$ is defined as follows, we set for $m\ge 2$
$$ sd^m (\sigma) = sd(sd^{m-1} (\sigma)).$$
$sd^0(\sigma)$ and $sd^1(\sigma)$ are defined respectively as $\sigma$ and $sd(\sigma)$.

\vspace*{0.3cm}

The following lemma is basic in Algebraic Topology.

\paragraph{{\sc Lemma 1.5.}}{\em
Given a simplex $\sigma$, and given $\epsilon > 0$, there is an $m$
such that each simplex of $ sd^m\sigma $ has diameter less than $\epsilon$.}

\vspace*{0.5cm}

For a proof of this lemma, see for example \cite{[Mu]}.
 
Let $\sigma =[\nu_1\ld \nu_p]$ and $\tau =[\mu_1\ld \mu_r]$.   We shall adopt the following notations
$$ [\sigma,\tau]=[\sigma,\mu_1\ld \mu_r]=[\nu_1\ld\nu_p,\tau]=[\nu_1\ld \nu_p,\mu_1\ld \mu_r].$$
Now let $\sigma $ be a $p$-simplex, $p\ge 2$, set
$$T(\sigma)= \Big[b(\sigma),\sigma\Big]+\sum^{p-2}_{\ell=1}\sum_{j_1 \ld j_{\ell}\atop {1\le j_i\le p-i+1}}   { (-1)^{j_1+\cdots +j_{\ell}}\Big[b(\sigma),b(\pl_{j_1}\sigma) \ld b(\pl ^{\ell}_{j_{\ell}\cdots j_1}\sigma) , \pl^{\ell}_{j_{\ell}\cdots j_1}\sigma\Big ]}$$
and extend $T$ by linearity to $S_p$.

If $\tau$ is an element of $S_1$ then we set 
$$T(\tau)=0$$
\paragraph{{\sc Proposition 1.6 }}{\em
If $\sigma$ is an element of $S_p$, $p\ge 2$, then
$$\pl {T(\sigma)} + T(\pl \sigma) = sd(\sigma) - \sigma .$$}

\vspace*{0.3cm}

\noindent This proposition follows by a straightforward computation.

\section{\hspace*{-1.5ex}Fundamental solutions for $\olp \sb {\rm b}$}

In this section, we shall construct two fundamental solutions for the tangential Cauchy-Riemann Complex. The  second solution  will be derived  from the first and will  yield optimal H\"older estimates for $\olp \sb {\rm b}$. 

Let us begin by some notations.

\paragraph{2.0. Notations.}
Throughout this section  $M$ will denote a  
$q$-concave $CR$ generic ${\cal C}^2$ submanifold of codimension $k$ in $\C^n$.

${\cal I}$ is the set of all subsets $I\subseteq\{\pm
1\ld\pm k\}$ such that $|i|\ne |j|$ for all $i,j\in I$ with $i\ne j$.

For $I\in{\cal I}$, $|I|$ denotes the number of elements in $I$ . We set 
$$ \Delta_{1\cdots |I|}=\{(\lambda_1,\ld,\lambda_{|I|})\in (\R^+)^{|I|} {~with~}  {\sum^ {|I|}_{j= 1}\lambda_j}=1\}$$
${\cal I}(\ell), 1\le\ell\le k$, is the set of all $I\in{\cal I}$ with
$|I|=\ell$.

\noindent ${\cal I}'(\ell), 1\le\ell\le k$, is the set of all $I\in{\cal I}(\ell)$
of  the form $I=\{j_1\ld j_\ell\}$ with $|j_\nu |=\nu $ for $\nu
=1\ld\ell$.

If $I\in{\cal I}$ and $\nu \in\{1\ld |I|\}$, then $I_\nu $ is the element with number $\nu $ in $I$  after ordering $I$ by modulus.We set  $I(\hat\nu)  =I\setminus\{I_\nu \}$ .    

If $I\in{\cal I}$, then
$$\sgn I :=\cases{
1 & if the number of negative elements in $I$ is even\cr
-1 & if the number of negative elements in $I$ is odd\cr}$$

\paragraph{2.1. First fundamental solution for $\olp \sb {\rm b}.$}
In this section we shall construct our first fundamental solution for the tangential Cauchy-Riemann complex.

\noindent Let  $z^0\in M$, 
$U\subseteq\C^n$ be a neighborhood of $z^0$ and  $\hat\rho _1\ld\hat\rho _k :
U\to\R$ be  functions of class ${\cal C}^2$ such that ~:
$$M\cap U=\{\hat\rho _1=\cdots =\hat\rho _k=0\}\hbox{~and~}\pl\hat\rho
_1(z^0)\ww\pl\hat\rho _k(z^0)\ne 0.$$ 
Since $M$ is $q$-concave , it follows from lemma 3.1.1 in
\cite{[Ai/He1]} that  we can find a constant $C>0$ such that the functions
$$\begin{array}{rl}
\rho _j &:= \hat\rho _j + C \sum^k_{\nu =1} \hat\rho ^2_\nu \quad
(j=1\ld k)\\
\rho _j &:= -\hat\rho _{-j} + C \sum^k_{\nu =1} \hat\rho ^2_\nu \quad
(j=-1\ld -k)\end{array}$$
have the following property~: for each  $I\in{\cal I}$  and every $\lambda
\in\Delta _{1\cdots |I|}$ the Levi form of $\lambda _1\rho
_{I_1}+\cdots +\lambda _{|I|}\rho _{I_{|I|}}$ at $z^0$  has at least $q+k$
positive eigenvalues.

\noindent Let $(e_1\ld e_k)$ be the canonical basis of $\R^k$, set
$e_{-j}:=-e_j$ for every $1\le j\le k$.

\noindent Let $I=(j_1\ld j_k)$ be in ${\cal I}'(k)$ ,  set 
$$ \Delta_{ I}=\{\sum^k_{i=1} \lambda_i e_{j_i} {~with~} \lambda_i\ge 0, {~all~} i, {~and~}\sum^k_{i=1} \lambda_i=1\},$$
and for each $a=\sum^k_{i=1} \lambda_i e_{j_i}$,  let $G_a$ and $\Phi_a$ be respectively the Leray map and the barrier function at $z^0$ corresponding to $\rho_a=\lambda_1\rho_{j_1}+\cdots+\lambda_k\rho_{j_k}$ (see sect. 1.4).
 
\noindent We call $\rho_a$ (resp. $\phi_a$) the defining function
(resp. the barrier function) of $M$ in the direction $a$.
 
\noindent Let $\sigma=[a^1\ld a^p]$, $p\ge 1$, be a collection of $p$ vectors , where $a^i\in {\bigcup_{I\in {\cal I}'(k)}}\Delta_I $, for every $1\le i\le k$.

\noindent Then we define 
$$\tilde\Omega[\sigma]:=\Omega(G_{a^1}\ld G_{a^p})$$
(cf sect 1.3) , and for every $0\le s\le n$ and every  $0\le r\le n-p$, we define $\tilde\Omega_{s,r}[\sigma]$ as the piece of $\tilde\Omega[\sigma]$ which is of type $(s,r)$ in $z$.

\noindent If we denote by $S'_p$ the set of all finite formal linear combinations of such collections, with integer coefficients,  and we extend $\tilde\Omega$ by linearity to $S'_p$; then  the generalized Koppelman lemma implies 
\paragraph{{\sc Lemma 2.1}}{\em 
For  every $\tau\in S'_p$, we have
$$\olp_{\zeta,z} \tilde\Omega[\tau]=\tilde\Omega [\pl\tau]$$
outside the singularities.}

\vspace*{0.5cm}

\noindent  Let $I=(j_1\ld j_l)$ be in ${\cal I}'(l)$,  $1\le l\le k$ and $\sigma_I=[e_{j_1}\ld e_{j_l}]$. Then by continuity of the Levi form, by lemma 1.3 and lemma 1.5, we can find a positive integer $m$ independant of $I$ and $l$ such that  for every simplexe $\tau =[a^1\ld a^l]$ in $sd^m(\sigma_I)$, the Leray maps of 
$G_{a^1}\ld G_{a^l}$ are $q+k$-holomorphic in the same directions with
respect to the variable $z\in \C^n$. Therefore we have the following lemma. 
\paragraph{{\sc Lemma 2.2}}{\em 
 There is a positive integer m such that for every $I\in {\cal I}'(l)$, $1\le l\le k$, any $s\ge 0$  and every $r\ge n-k-q+1$
 
\decale {(i)} $\tilde\Omega_{s,r}(sd^m (\sigma_I))=0$

\decale {(ii)}$\olp_z \tilde\Omega_{s,r-1}(sd^m (\sigma_I))=0$

\noindent on the set where all the denominators are non-zero.}

\vspace*{0.5cm}

{\sl Proof. } this follows by linearity from the fact that
$\Omega(G_{a^1}\ld G_{a^l})=0$, for any $[a^1\ld a^l]$ in $sd^m(\sigma_I)$. The last statement is easy to prove , looking at the definition of $\Omega$ (see sect. 1.3), because 
$G_{a^1}\ld G_{a^l}$ are $q+k$-holomorphic in the same directions with
respect to the variable $z\in \C^n$

\noindent By the same arguments, we have
$\tilde\Omega_{s,r}(\partial(sd^m
(\sigma_I)))=\tilde\Omega_{s,r}(sd^m(\partial\sigma_I)=0$, for all
$r\ge n-k-q+1$, and from lemma 2.1, we have
$$\overline\partial_z\tilde\Omega_{s,r-1}(sd^m(\sigma_I))=-\overline\partial_\zeta\tilde\Omega_{s,r}(sd^m(\sigma_I))+\tilde\Omega_{s,r}(\partial(sd^m(\sigma_I))$$ which implies the statement (ii).\hfill$\Box$

\vspace*{0.5cm}

Now let $D$ be a neighborhood of $z^0$ such that for every
$1\le l\le k$ , all $0\le i\le m$ and every vertex $a$ in $sd^i(\sigma_I)$ , the barrier function $\Phi_a$ satisfies an inequality such  $(1.5)$ , for $\zeta,z \in D$. Set

 $$M_0 := M\cap D, $$
and for $I\in{\cal I}$
$$\begin{array}{rl}
D_I &:=\{\rho _{I_1}<0\}\cap\cdots\cap\{\rho _{I_{|I|}}<0\}\cap D,\\
D^*_I &:=\{\rho _{I_1}>0\}\cap\cdots\cap\{\rho _{I_{|I|}}>0\}\cap D,\\
S_I &:=\{\rho _{I_1}=\cdots =\rho _{I_{|I|}}=0\}\cap D,\\
S^+_{\{j\}} &:=\overline {D^*}_{\{j\}}\hbox{~for~} j=\pm 1\ld \pm k,\\
S^+_I &:= S_{I(|\widehat{I}|)}\cap\overline {D^*}_{\{I_{|I|}\}}\hbox{~if~} I\in{\cal I}\hbox{~and~}|I|\ge 2.\end{array}$$
We oriente these manifolds as follows~:
$$\begin{array}{rl}
D_I \hbox{~and~} D^*_I &\hbox{~as~}\C^n\quad\forall I\in{\cal I}\\
S^+_{\{j\}} &\hbox{~as~} D^*_{\{j\}}\hbox{~for~} j=\pm 1\ld\pm k\\
S_I &\hbox{~as~} \pl S^+_I\quad I\in{\cal I}\\
S^+_I &\hbox{~as~} S_{I(|\widehat{I}|)} \hbox{~for all~} I\in{\cal I}\hbox{~such that~}
|I|\ge 2\\
M_0 &\hbox{~as~} S_I\hbox{~where~} I=\{1\ld k\}.\end{array}$$
Fix $1\le l\le k$ and $I\in {\cal I}'(l)$.

\noindent Let $B=\big(\overline{\zeta_1-z_1}\ld\overline{\zeta_n-z_n}\big)$ and define
$$\tilde\Omega_B[\tau]:=\Omega(B,G_{\nu^1}\ld G_{\nu^p})\eqno (2.1)$$ for any  $\tau=[\nu^1\ld \nu^p]$ in $S'_p$, $p\ge 1$.   Extend this operation, by linearity,  to all elements of $S'_p$.  

\noindent Now by applying  lemma 1.1 we get
$$\olp_{\zeta,z}\tilde\Omega_B[\sigma_I]=-\tilde\Omega[\sigma_I]-\tilde\Omega_B[\partial\sigma_I]\eqno (2.2)$$
(where $\tilde\Omega_B[\partial\sigma_I]:=\Omega(B)$ if $|I|=1$) for $z\in\overline {D_I}$ and $\zeta\in\overline {D^*_I}$, with $\zeta\ne z$.

\noindent Let $|I|\ge 2$ and  $T$ be defined as in sect. 1.5 and $m$ an integer such that lemma 2.2 holds.

\noindent By applying lemma 1.4, lemma 2.1 , proposition 1.6 , we obtain
$$\olp_{\zeta,z}\sum^{m-1}_{i=0}\tilde \Omega[T(sd^i(\sigma_I))]=-\tilde\Omega[\sigma_I]- \sum^{m-1}_{i=0}\tilde \Omega[T(sd^i(\pl\sigma_I))]+\tilde \Omega[sd^m(\sigma_I)]\eqno(2.3)$$
for $z\in\overline {D_I}$ and $\zeta\in\overline {D^*_I}$, with $\zeta\ne z$.

\noindent Now define for $|I|\ge 2$
$$ K^I(\zeta,z)=\tilde\Omega_B[\sigma_I](\zeta,z)-\sum^{m-1}_{i=0}\tilde \Omega[T(sd^i(\sigma_I))](\zeta,z),$$
$$B^I(\zeta,z)=\sum^{|I|}_{\nu=1}(-1)^{\nu+1} K^{I(\hat\nu)}=-\tilde\Omega_B[\partial\sigma_I](\zeta,z)+\sum^{m-1}_{i=0}\tilde \Omega[T(sd^i(\partial\sigma_I))](\zeta,z)\eqno (2.4)$$
and set for $|I|=1$
$$K^I(\zeta,z)=\tilde\Omega_B[\sigma_I](\zeta,z),$$and
$$B^I(\zeta,z)=\Omega(B)(\zeta,z)\hbox {~(the~ Martinelli-Bochner- Koppelman~ kernel).~}\eqno (2.5)$$
Then we have the following

\paragraph{{\sc lemma 2.3}}{\em  

\vspace*{0.3cm}

\decale{(i)} For $z\in\overline {D_I}$ and $\zeta\in\overline {D^*_I}$, with $\zeta\ne z$, we have
$$\olp_{\zeta,z}K^I=B^I-\tilde \Omega[sd^m(\sigma_I)]$$
\decale{(ii)} There exist a constant ${\cal C}>0$ and a finite family $\{\gamma_1\ld\gamma_L \}$ of linearly independent families $\gamma_i=[\gamma^1_i\ld\gamma^{|I|}_i ]$ in $\Delta_I$ such that
 
$$\|K^I(\zeta,z)\|\le{\cal C }\sum^L_{i=1} {1\over \prod^{|I|}_{j =
    1}|\Phi^{\gamma^j_i}(\zeta,z)|{~} | \zeta-z|^{2n-2|I|-1}}$$}

\vspace*{0.5cm}

{\sl Proof.} (i) is a consequence of (2.1), (2.3) and (2.4). The
estimate in (ii) is easy to see from the definition of $\Omega$, by
using lemma 1.2 and inequality (1.5) (cf. the proof of Lemma
2.6).\hfill$\Box$

\vspace*{0.5cm}

The following lemma shows that the kernel $K^I$ (resp. $B^I$) has
locally integrable coefficients on $S_I$ (resp. $S^+_I$) in both variables $\zeta$ and $z$.

\paragraph{{\sc  Lemma 2.4}}{\em 

\vspace*{0.3cm}

\decale{(i)} Let $I\in{\cal I}$ and $(\gamma ^1\ld\gamma ^{|I|})$ be a
family of linearly independent vectors in $\R^{|I|}$ and $z\in\old_I$
then there exists $C>0$ and $\varepsilon _0>0$ such that for all
$\varepsilon <\varepsilon _0$ and for all $j\in\{\pm 1\ld \pm k\}\setminus I, |j|\le |I|$

$${\displaystyle \int_{{\zeta\in S_I}\atop{|\zeta -z|<\varepsilon }} { d\lambda (\zeta)\over
\prod^{|I|}_{i=1}|\Phi _{\gamma^i}(\zeta,z)||\zeta -z|^{2n-2|I|-1}}  \le C \varepsilon
\big(1+|\ln\varepsilon |\big)^{|I|}} \eqno (2.6)$$

$${\displaystyle\int_{{\zeta\in S^+_{I\cup\{j\}}}\atop{|\zeta -z|<\varepsilon }} { d\lambda (\zeta)\over
\prod^{|I|}_{i=1}|\Phi _{\gamma
^i}(\zeta,z)||\zeta -z|^{2n-2|I|-1}} \le C \varepsilon
\big(1+|\ln\varepsilon |\big)^{|I|} }\eqno (2.7)$$

$${\displaystyle\int_{{\zeta\in S^+_I}\atop{|\zeta -z|<\varepsilon }} { d\lambda (\zeta)\over
\prod^{|I|-1}_{i=1}|\Phi_{\gamma^i}(\zeta,z)||\zeta -z|^{2n-2|I|+1}}  \le C \varepsilon
\big(1+|\ln\varepsilon |\big)^{|I|-1}}\eqno  (2.8)$$

\decale{(ii)}All of the above estimates hold if we integrate with respect to $z$ instead of $\zeta$.}

\vspace*{0.5cm}
 
{\sl Proof.} 
 Since $M$ is $CR$ generic and $\gamma ^1\ld\gamma ^{|I|}$ are
 linearly independent, we can take $Im\Phi_{\gamma^1}(\cdot,z)\ld
 Im\Phi_{\gamma^|I|}(\cdot,z)$ as  coordinates on $S_I$ and
 $S^+_{I\cup\{j\}}$, for $|z-\zeta|<\varepsilon $ with $\varepsilon>0$
 sufficiently  small (cf. lemma 2.3 in \cite{[Ba2]}). Thus taking into account the following inequality (cf. (1.5))
$$|\Phi_{\gamma^i}(\zeta,z)|\ge C
(|Im\Phi_{\gamma^i}(\zeta,z)|+|\zeta-z|^2)$$  for $\zeta, z\in M_0$
with $|\zeta-z|<\varepsilon$, we see that the left hand side term in  $(2.6)$ (resp.$(2.7)$) is bounded by
$$
\int_{{X\in \R^{2n-|I|}}\atop{|X|<\varepsilon }} { dX\over\prod^{|I|}_{i=1}\big(|X_i|+|X|^2\big)|X|^{2n-2|I|-1}} \le C \varepsilon  \big(1+|\ln\varepsilon |\big)^{|I|}.$$
$(2.8)$ is proved likewise and the proof of (ii) is similar.\hfill$\Box$

\paragraph{{\sc  Definition 2.5}}{\em

\vspace*{0.3cm}

\decale{(i)} Let $I\in{\cal I}$ and $j\in\{\pm 1\ld\pm k\}\setminus I, |j|\le|I|$.
It follows from lemma 2.3  and estimates $(2.6)$ , $(2.7)$, that the following operators are well defined and continuous~:
{\parindent=0pt
$$\begin{array}{rll}
\widehat{K}^I_{0,r} &: {\cal C}^0_*(S_I)  \longrightarrow
{\cal C}^0_{0,r}(\old_I)\cap{\cal C}^\infty_*(D_I),&\\
\widehat{K}^{I,j}_{0,r} &: {\cal C}^0_*(S^+_{I\cup\{j\}})  \longrightarrow
{\cal C}^0_{0,r}(\old_I)\cap{\cal C}^\infty_*(D_I)&\\

\noalign{where}\\

\widehat{K}^I_{0,r}f(z) &:=\int_{\zeta\in S_I} f(\zeta)\wedge K^I_{0,r}(\zeta,z),\quad
z\in\old_I, f\in{\cal C}^0_*(S_I),&\\

\widehat{K}^{I,j}_{0,r}f(z) &:=\int_{\zeta\in S^+_{I\cup\{j\}}}
f(\zeta)\wedge K^I_{0,r}(\zeta,z),\quad z\in\old_I, f\in{\cal C}^0_*(S^+_{I\cup\{j\}})&\end{array}$$}
\decale{(ii)} Let $I\in{\cal I}$, it follows from lemma 2.3 and estimate $(2.8)$  that the operator defined by~:
$$\widehat B^I_{0,r} f(z):=\int_{\zeta\in S^+_I}f(\zeta)\wedge B^I_{0,r}(\zeta,z),\quad z\in\old_I, f\in{\cal C}^0_*(S^+_I)$$
is continuous from ${\cal C}^0_*(S^+_I)$ into
${\cal C}^0_*(\old_I)\cap{\cal C}^\infty_*(D_I)$.}

\paragraph{{\sc Remark.}} {\em If $|I|\ge 2$ et $\nu \in\{1\ld |I|\}$ then  $S^+_{I(\hat\nu )\cup(I_\nu )}=S^+_I$ and therefore
$$\widehat{K}^{I(\hat\nu ),I_\nu } f(z)=\int_{\zeta\in S^+_I} f(\zeta)\wedge K^{I(\hat\nu)}_{0,r}(z,\zeta)\hbox{~~for~~} z\in\old_{I(\hat\nu )}\hbox{~and~}
f\in{\cal C}^0_*(S^+_I)$$
thus from  the  definition of $B^I_{0,r}$, we obtain
$$\widehat B^I_{0,r} f(z)=\sum^{|I|}_{\nu =1}(-1)^{\nu +|I|}\widehat{K}^{I(\hat\nu ),I_\nu}_{0,r}f(z).\eqno (2.9)$$}
 
\paragraph{{\sc Lemma 2.6}}{\em Let  $I\in{\cal I}, n-k-q+1\le r\le n-k$
and $f\in {\cal C}^1_{0,r}(D)$ with compact support on $D$. Then  the following
equality holds  in the sens of currents~:
$$\olp\widehat{K}^I_{0,r-1} f+(-1)^{|I|+1}\widehat{K}^I_{0,r}
\olp f=(-1)^{r}\big[\olp\widehat{B}^I_{0,r-1}f+(-1)^{|I|}\widehat{B}^I_{0,r}\olp f\big]\eqno(2.10)$$
on  $\old_I$.}
 
\vspace*{0.5cm}
 
{\sl Proof.} The following identity is true from lemma 2.3(i) and Stokes' theorem~:
if  $n-k-q+1\le s\le n-k, g\in C^1_{0,s+1}(D)$ with compact support on $D$ and
if $z\in D_I$ then~:
$$\begin{array}{rl}
\widehat{K}^I_{0,s}g(z)&=\int_{S^+_I}\olp g\wedge
K^I_{0,s}(\cdot,z)+(-1)^{s+1}\widehat{B}^I_{0,s}g(z) + (-1)^{|I|}\\
 &\olp \int_{S^+_I}g\wedge K^I_{0,s-1}(\cdot,z)+(-1)^{s}\int_{S^+_I} g\wedge \tilde\Omega_{0,s}(sd^m{\sigma_I})(\cdot,z).\end{array}\eqno(2.11)$$
Since the forms
$\widehat{K}_{0,r-1}f,\widehat{K}_{0,r-1}\olp f,\widehat{B}_{0,r-1}f$ and $\widehat{B}_{0,r}\olp f$ are
continuous on $\old_I$, it is sufficient to prove  $(2.10)$ on $D_I$
where  these forms are smooth.

\vspace*{0.3cm} 

By setting $s=r$ and $g=\olp f$ in $(2.11)$ and using lemma 2.2 (i), we obtain
$$\widehat {K}^I_{0,r} \olp f(z)=(-1)^{r+1}\widehat{B}^I_{0,r} \olp f(z)+(-1)^{|I|}\olp \int_{S^+_I}\olp f\wedge
K^I_{0,r-1}(\cdot,z)\eqno (2.12)$$
for all $z\in D_I$.
 
\noindent If we set now $s=r-1$  and  $g=f$ in $(2.11)$, then we get
\begin{eqnarray*}
\widehat{K}^I_{0,r-1}f(z)&=&\int_{S^+_I}\olp f\wedge
K^I_{0,r-1}(\cdot,z)+(-1)^{r}\widehat{B}^I_{0,r-1}f(z)\\
&&+(-1)^{|I|}\olp\int_{S^+_I}f\wedge K^I_{0,r-2}(\cdot,z)+(-1)^{r-1}\int_{S^+_I} f\wedge \tilde\Omega_{0,r-1}(sd^m{\sigma_I})(\cdot,z),
\end{eqnarray*}
and then by lemma 2.2 (ii) 
$$\olp\widehat{K}^I_{0,r-1}f(z)=\olp\int_{S^+_I}\olp f\wedge K^I_{0,r-1}(\cdot,z)+(-1)^{r}
\olp\widehat{B}^I_{0,r-1}f(z)\eqno (2.13)$$
for all $z\in D_I$.
 
\noindent The lemma now follows from $(2.12)$ and $(2.13)$ .\hfill$\Box$

\vspace*{0.5cm} 

Now define
$$K(\zeta,z):=\sum_{I\in {\cal I}'(k)} (sgnI) K^I(\zeta,z)\eqno (2.14)$$
for $\zeta,z\in M_0$ with $\zeta\ne z$, and denote by $K_{s,r}$ the
piece of $K$ which is of type $(s,r)$ in  the variable $z$.  

\noindent From lemma 2.4, we see that the kernel $K$ has locally integrable coefficients in both variables $\zeta$ and $z$.

\noindent Now by applying (2.10) $k$ times , taking into account (2.9) and using the classical Martinelli-Bochner-Koppelman formula  (see \cite{[Ba1]} or  \cite{[Ba2]} for technical details) we obtain  the following integral representation
\paragraph{{\sc Theorem 2.7.}}{\em

Let $\Omega\sub M_0$ of piecewise ${\cal C}^1$ boundary and $f$ a
$(0,r)$  ${\cal C}^1$ form on $\ol\Omega$ with  $n-k-q+1\le r\le n-k$ , then
$$\displaylines{
(-1)^{r(k+1)}f(z)=\int_{b\Omega }f(\zeta)\wedge K_{0,r}(\zeta,z)-\int_\Omega \overline\partial\sb {\rm b}f(\zeta)\wedge K_{0,r} (\zeta,z)\hfill\cr 
\hfill +(-1)^{k+1}  \overline\partial\sb {\rm b} \int_{\Omega }f(\zeta)\wedge K_{0,r-1}(\zeta,z).\cr}$$}
 By a duality argument we obtain
\paragraph{{\sc Corollary 2.8}}{\em
 
Let  $\Omega \sub M_0$ of piecewise   ${\cal C}^1$ boundary and 
   $f$ a ${\cal C}^1$  $(0,r)$-form on $\ol\Omega$
with  $0\le r\le q-1$, then  we have  
$$\displaylines{
(-1)^{ r(k+1)} f(\zeta) = \int_{z\in b\Omega } f(z)\wedge
K_{n,n-k-1-r}(\zeta,z)-\int_{z\in\Omega } \overline\partial\sb {\rm b} f(z)\wedge K_{n,n-k-1-r}(\zeta,z)+\hfill\cr
\hfill +(-1)^{(k+1)} \overline\partial\sb {\rm b} \int_\Omega f(z)\wedge K_{n,n-k-r}(\zeta,z).,\cr}$$}

\noindent We say that $K$ is a fundamental solution for $\overline\partial\sb{\rm b}$ on $M_0$.

\paragraph{2.2. Second  fundamental solution for $\olp \sb {\rm b}.$}
In this section, we shall construct our second fundamental solution for the tangential Cauchy-Riemann complex on $M_0$.
This fundamental solution will be derived from the first one , by
using an idea of Henkin \cite {[He3]} (cf. \cite {[Bo]}).

\vspace*{0.3cm}

Let $m$ be as in Lemma 2.2 and $\nu^*\in\bigcup_{I\in {\cal I}'(k)}\Delta_I $ such that

$$ \cases
{\vbox{\hsize=10cm{\noindent For any  $k$- simplex $\tau$ in $sd^m(\sigma_I)$,
each  collection  of $k$ elements in $[\nu^*,\tau]$ is a $k$-simplex. 
  }}\cr}\leqno(*)$$
\paragraph{{\sc Remark. }}{\em
The choice of such $\nu^*$ is very important for our optimal estimates.}  

\vspace*{0.3cm}

\noindent We adopt the following notation
$$ [\nu^*, \sum_i c_i\sigma_i]=\sum_ic_i[\nu^*,\sigma_i]$$
for any element $\sum_i c_i\sigma_i$ in $S'_p$.

\noindent Set
$$E(\zeta,z)=\sum_{I\in{\cal I}'(k)}sgnI\Big( \tilde\Omega_B[\nu^*,\sigma_I]+ \sum^{m-1}_{i=0}\tilde \Omega[\nu^*,T(sd^i(\sigma_I))]\Big)\eqno (2.15)$$
and
$$R(\zeta,z)=\sum_{ I\in {\cal I}'(k) }( sgnI)\tilde\Omega[\nu^*,sd^m(\sigma_I)].\eqno (2.16)$$
Since
$$\sum_{I\in{\cal I}'(k)}(sgnI)\pl\sigma_I=0, $$
then by applying lemma 1.1, proposition 1.6 and  (2.2), we obtain 
$$\olp_{\zeta,z}E(\zeta,z)=K(\zeta,z)-R(\zeta,z)\eqno (2.17)$$
for $\zeta,z\in M_0$ with $\zeta\ne z.$

\noindent Now we claim that $R$ is a fundamental solution for
$\overline\partial\sb {\rm b}$ on $M_0$, this means that Theorem 2.7
holds also for the kernel $R$. To prove it, following  Henkin \cite
{[He3]}, all we have to do is to show that the singularity of E is
mild enough so that the identity (2.17) holds on all $M_0 \times M_0$
in the sens of distributions. For once this is done, our claim follows
by applying $\olp_{\zeta,z}$ to both sides of (2.17) and then using Theorem 2.7.

\noindent The proof of the first part of Theorem 0.1 will be then  complete by setting
\paragraph{{\sc Definition 2.9}}
$$ {\cal R}_r(\zeta,z) :=\cases{
{(-1)}^{r(k+1)}R_{0,r}(\zeta,z) & if $n-k-q\le r\le n-k$ \cr
{(-1)}^{r(k+1)}R_{n,n-k-1-r}(z,\zeta) & if $0\le r\le q-1$ .\cr}$$
 
\vspace*{0.5cm}

Now to realize our program, we follow the proof of Theorem 1, chap.21  in \cite{[Bo]}.

\noindent First we need the following lemma
\paragraph{{\sc Lemma 2.10.}}{\em
Given $W\sub M_0$, there is a positive constant $C$ such that for each $\epsilon >0$ and $z\in W$, we have

\decale {(i)}$\displaystyle{\int_{\zeta\in M_0\atop |\zeta-z|\le \epsilon}\|K(\zeta,z)\|{~}d\lambda(\zeta)}\le C\epsilon (1+|\ln\epsilon|)^k$ 

\decale{(ii)}$\displaystyle{\int_{\zeta\in M_0\atop |\zeta-z|\le \epsilon}\|R(\zeta,z)\|{~}d\lambda(\zeta)}\le C{\epsilon}$

\decale{(iii)}$\displaystyle{\int_{\zeta\in M_0\atop |\zeta-z|\le
    \epsilon}\|E(\zeta,z)\|{~}d\lambda(\zeta)}\le C{\epsilon^2}
(1+|\ln\epsilon|)^k.$

\decale{$(i\nu)$} All of the above inequalities hold if we integrate with respect to $z$ instead of $\zeta$.}

\vspace*{0.5cm}

\noindent Let us assume the lemma for the moment  and show that equation (2.17) holds on all $M_0\times M_0$. 

\noindent For $\epsilon>0$, choose a smooth function $\chi_\epsilon$ on  $M_0\times M_0$ with the following properties~
$$\chi_\epsilon(\zeta,z)   =\cases{
1 & if $|\zeta-z|\ge\epsilon$\cr
0 & if $|\zeta-z|\le{\epsilon\over 2}$.\cr}$$
and for any first-order derivative $\cal D$,
$$|\cal D\{\chi_\epsilon\}|\le {C\over\epsilon}\eqno (2.18)$$
where $\cal C$ is a positive constant that is independent of $\epsilon$.

\noindent Since $\chi_\epsilon$ vanishes near the diagonal of $M_0\times M_0$, we have from (2.17)
$$\olp_{\zeta,z}\{\chi_\epsilon E\}=(\olp_{\zeta,z}\chi_\epsilon)\wedge E + \chi_\epsilon (K - R)\eqno (2.19)$$
on $M_0\times M_0$. From Lemma 2.10, we have

$$\chi_\epsilon K\rightarrow K, \chi_\epsilon R\rightarrow R
, \chi_\epsilon E\rightarrow E \hbox {~and~} \olp_{\zeta,z}\{\chi_\epsilon E\}\rightarrow \olp_{\zeta,z}E $$
in the sens of currents, as $\epsilon\rightarrow 0$ 

\noindent From part (iii) in lemma 2.10  and  estimate (2.18), we see that
$$(\olp_{\zeta,z}\chi_\epsilon)\wedge E \rightarrow 0$$
as $\epsilon\rightarrow 0$, in the sens of currents.
So we obtain the desired result by letting  $\epsilon\rightarrow 0$ in the equation (2.19).
   
\paragraph{{\sc Proof of lemma 2.10.}}

Looking at the definitions of the kernels $K$, $E$ and $R$ (cf.
(2.14),...,(2.16)) and taking into account Lemma 1.2, we see that we have to estimate the following typical term

$${{\cal N}(\zeta,z)\over \prod^k_{i = 1}(\Phi_{a^i}(\zeta,z))^{r_i}(\Phi_{a^0}(\zeta,z))^{r_0}(\Phi_{a^{k+1}}(\zeta,z))^{r_{k+1}}{~} |\zeta-z|^{2s}}\eqno (2.20)$$
where $a^1\ld a^k$ are linearly independent, $a^0=\sum^{k}_{i=1}x_i a^i$, $a^{k+1}=\sum^{k}_{i=1}y_i a^i$,  
$$r_i\ge 1, \hbox{~all~}1\le i\le k;{~~}  s,r_0, r_{k+1}\ge 0\hbox{~and~}$$
$$s+\sum^{k+1}_{i=0}r_i=n.$$

For the kernel $K$, we have $r_{k+1}=0$ and  

\noindent either $r_0=0$, $s\ge 1$ and the function $\cal N$ involves coefficients of the differential form
$$\Big(G_{a^1}.d(\zeta-z)\Big)\wedge\cdots\wedge\Big(G_{a^k}.d(\zeta-z)\Big)\wedge\Big(\overline{(\zeta-z)}.d(\zeta-z)\Big)$$
or $s=0$ , $r_0\ge 1$ and  the function $\cal N$ contains the coefficients of the term
$$\Big(G_{a^1}.d(\zeta-z)\Big)\wedge\cdots\wedge\Big(G_{a^k}.d(\zeta-z)\Big)\wedge\Big(G_{a^0}.d(\zeta-z)\Big)$$
Since
$$G_{a^0}(\zeta,z)=\sum^{k}_{i=1}x_iG_{a^i}+{\cal O}(|\zeta-z|)$$
we obtain in both cases
$$|{\cal N}(\zeta,z)|\le C|\zeta-z|\eqno (2.21)$$
Since $M$ is $CR$ generic and $a^1\ld a^k$ are linearly independent,

\noindent $Im\Phi_{a^1}(.,z)\ld Im\Phi_{a^k}(.,z)$ can be taken as
local coordinates on $M_0$ (cf.Lemma 2.3 in \cite{[Ba1]}). Then in view of inequality (1.5), 
$$\begin{array}{rlr}\int_{\zeta\in M_0\atop |\zeta-z|\le
    \epsilon}\|K(\zeta,z)\|{~}d\lambda(\zeta)&\le C \int_{X\in\R^{2n-k}\atop{|X|<\epsilon}}
{dX\over\prod^k_{j=1}\big(|X_j|+|X|^2\big)|X|^{2n-2k-1}}&\\

&\le C\epsilon(1+|\ln\epsilon|)^k.&\end{array}$$

Now for the kernel $E$, we have $r_{k+1}\ge 1$ and

\noindent either $s=0$, $r_0\ge 1$ and the function $\cal N$ involves the coefficients of the term
$$\Big(G_{a^1}.d(\zeta-z)\Big)\wedge\cdots\wedge\Big(G_{a^k}.d(\zeta-z)\Big)\wedge\Big(G_{a^0}.d(\zeta-z)\Big)\wedge\Big(G_{a^{k+1}}.d(\zeta-z)\Big)$$
\noindent or $s\ge 1$, $r_0=0$ and the the function $\cal N$ contains the coefficients of the differential form
$$\Big(G_{a^1}.d(\zeta-z)\Big)\wedge\cdots\wedge\Big(G_{a^k}.d(\zeta-z)\Big)\wedge\Big(\overline{(\zeta-z)}.d(\zeta-z)\Big)\wedge\Big(G_{a^{k+1}}.d(\zeta-z)\Big)$$
By the same arguments as above, we obtain in this case
$$|{\cal N}(\zeta,z)|\le C|\zeta-z|^2$$
and
$$\begin{array}{rlr}\int_{\zeta\in M_0\atop |\zeta-z|\le \epsilon}\|E(\zeta,z)\|{~}d\lambda(\zeta)&\le \int_{X\in\R^{2n-k}\atop{|X|<\epsilon}}
{dX\over\prod^k_{j=1}\big(|X_j|+|X|^2\big)|X|^{2n-2k-2}}&\\

&\le C {\epsilon^2}(1+|\ln\epsilon|)^k.&\end{array}$$

For the kernel $R$, we have $s=0$, $r_0=0$,  $r_{k+1}\ge 1$ and every
collection of $k$ elements in $\{a^1\ldots a^{k+1}\}$ is a family of
linearly independent vectors (see condition $(*)$ and the remark just
below  ).

\noindent  First it is easy to see , just as above, that  inequality (2.21) holds also in this case.

\noindent On other hand the following inequality is easy to prove: if $0\le\alpha_1\ld\alpha_k\le\alpha_{k+1}$, then
$$\prod^{k+1}_{i=1}\alpha_i\ge\prod^k_{i=1}\alpha^{1+{1\over
    k}}_i\eqno (2.22)$$
If we use (2.22) with $\alpha_i=|\Phi_{a^i}|$, $1\le i\le k+1$ ( up to
a permutation of $\{1\ld{k+1}\}$), then  by using local coordinates as
    above and inequality(1.5), we obtain
$$\begin{array}{rlr}\int_{\zeta\in M_0\atop |\zeta-z|\le
    \epsilon}\|R(\zeta,z)\|{~}d\lambda(\zeta)&\le C \int_{X\in\R^{2n-k}\atop{|X|<\epsilon}}
{dX\over\prod^k_{j=1}\big(|X_j|+|X|^2\big)^{1+{1\over k}}|X|^{2n-2k-3}}&\\

&\le C \epsilon.&\end{array}$$
Thus the proof of $(i), (ii), (iii)$ in Lemma 2.10 is complete. $ (i\nu)$ follows in the same way.$\hfill\Box$

\section{\hspace*{-1.5ex}End of proof of Theorem 0.1}
In this section we shall prove ${\cal C}^{\ell+{1\over 2}}-$estimates.
We first prove  ${\cal C}^{1\over 2}-$estimates and then we derive
${\cal C}^{\ell+{1\over 2}}-$estimates by using a kind of integration
by parts argument (see  \cite{[Ma/Mi]} and \cite{[Fi]}).

\paragraph{3.1. ${\cal C}^{1\over 2}-$Estimates.}

Suppose $M$ of class ${\cal C}^2$. Recall from the previous section that the coefficients of the kernel $R(\zeta,z)$ have the form 
 
$${{\cal N}(\zeta,z)}\over \prod^{k+1}_{i = 1}(\Phi_{a^i}(\zeta,z))^{r_i}$$
where $a^1\ld a^{k+1}$ are vectors in $\R^k$ such that every subset of $k$ elements in $\{a^1\ld a^{k+1}\}$ is a family of linearly independent vectors (condition $(*)$), the estimate $(2.21)$ holds for $\cal N$ and 

$$r_i\ge 1, {~all~}1\le i\le k+1;{~}\sum^{k+1}_{i=0}r_i=n.$$
We have 
$$\int_{\zeta\in M_0}\|R(\zeta,z^1)-R(\zeta,z^2)\|{~}d\lambda(\zeta)\le J_1(z^1,z^2)+J_2(z^1,z^2)$$
where
$$J_1(z^1,z^2):=\int_{{\zeta\in M_0}\atop{|\zeta-z^1|\le |z^1-z^2|^{1\over 2}}}\big(\|R(\zeta,z^1)\|+\|R(\zeta,z^2)\|\big){~}d\lambda(\zeta)$$
and
$$J_2(z^1,z^2):=\int_{{\zeta\in M_0}\atop{|\zeta-z^1|\ge |z^1-z^2|^{1\over 2}}}\|R(\zeta,z^1)-R(\zeta,z^2)\|{~}d\lambda(\zeta)$$
It follows from lemma 2.10 (ii) that
$${J_1(z^1,z^2)}\le C |z^1-z^2|^{1\over 2}.$$
Since ${\cal N}(\zeta,z)$ is smooth in $z$, it is not difficult to see by the same arguments as in the proof of Lemma 2.10 that
$$\begin{array}{rlr}{J_2(z^1,z^2)} &\le C|z^1-z^2|\displaystyle\int_{{X\in\R^{2n-k}}\atop{|X|\ge |z^1-z^2|^{1\over 2 }}}{dX\over\big(|X_1|+|X|^2\big)^{2+{1\over k}}\prod^k_{j=2}\big(|X_j|+|X|^2\big)^{1+{1\over k}}|X|^{2n-2k-3}} &\\
&&\\
&\le C |z^1-z^2|^{1\over 2}.&\end{array}$$
Thus
$$\int_{\zeta\in M_0}\|R(\zeta,z^1)-R(\zeta,z^2)\|{~}d\lambda(\zeta)\le  C |z^1-z^2|^{1\over 2}.\eqno (2.23)$$
Analogously we can show that 
 $$\int_{z\in M_0}\|R(\zeta^1,z)-R(\zeta^2,z)\|{~}d\lambda(z)\le C|\zeta^1-\zeta^2|^{1\over 2}.\eqno (2.24)$$
under the hypothesis that $M$ is of class ${\cal C}^3$.
This is because $R(\zeta,z)$ involves second-order derivatives in $\zeta$ of
the defining functions of $M$.
\paragraph{3.2. ${\cal C}^{\ell+{1\over 2}}-$Estimates.}
We assume  that $M$ is of class ${\cal C}^{\ell+2}$ ($\ell\ge 1$).

\noindent Let $a^1\ld a^k$ be linearly independent vectors in ${\bigcup_{I\in {\cal I}'(k)}}\Delta_I $ and  $a^{k+1}=\sum^{k}_{i=1}y_i a^i$ with $y_i\ne 0$, all $1\le i\le k$ (this means that every collection of $k$ vectors in $\{a^1\ld a^{k+1}\}$ is a family of linearly independent vectors).

\noindent Denote by $\tilde\rho_i$ (resp. $\phi_i$) the defininig
function (resp. the barrier function) of $M$ in the direction $a^i$ for $1\le i\le k+1$.

\noindent ${\cal E}^j$($j\ge 0$)will denote a smooth differential form on $M\times M$ vanishing of order $j$ for $\zeta=z$.
It is clear that
$$\phi_{k+1}={\sum^{k}_{i=1}y_i \phi_i}+{\cal E}^2\eqno (2.25)$$
We need the following lemma.
\paragraph{{\sc Lemma 3.1}}{\em 
There exist $Y^{\zeta}_1\ld Y^{\zeta}_k$, tangential vector fields  to $M$ such that for every $\zeta\in M_0$ and every $1\le i,j\le k$,
$$Y^{\zeta}_i\phi_j(\zeta,\zeta)=\delta_{ij},$$
where $\delta {ij}$ is Kronecker's symbol.}

\vspace*{0.5cm}

{\sl Proof.} Since $M$ is $CR$ generic and $a^1\ld a^k$ are linearly independent, we have
$$\partial{\tilde\rho_1}\wedge\cdots\wedge\partial{\tilde\rho_k}\ne 0\mbox {~on~} M_0.$$
Then the matrix 
\[A=
\left(
\begin{array}{ccc}
<\partial{\tilde\rho_1}(\zeta),\partial{\tilde\rho_1}(\zeta)>&\cdots&<\partial{\tilde\rho_k}(\zeta),\partial{\tilde\rho_1}(\zeta)>\\
\vdots&&\vdots\\
<\partial{\tilde\rho_1}(\zeta),\partial{\tilde\rho_k}(\zeta)>&\cdots&<\partial{\tilde\rho_k}(\zeta),\partial{\tilde\rho_k}(\zeta)>

\end{array}
\right)
\]
is invertible for all $\zeta\in M_0$ ( here $<.,.>$ denotes the Hermitian inner product), and there exist $\nu_1\ld \nu_k\in\{1\ld n\}$ such that the matrix
\[B=
\left(
\begin{array}{ccc}
{\partial{\tilde\rho_1}\over\partial\overline\zeta_{\nu_1}}(\zeta)&\cdots&{\partial{\tilde\rho_k}\over\partial\overline\zeta_{\nu_1}}(\zeta)\\
\vdots&&\vdots\\
{\partial{\tilde\rho_1}\over\partial\overline\zeta_{\nu_k}}(\zeta)&\cdots&{\partial{\tilde\rho_k}\over\partial\overline\zeta_{\nu_k}}(\zeta)
 
\end{array}
\right)
\]
is also invertible for all $\zeta\in M_0$.

\noindent Set 
$$Y^{\zeta}_i={1\over 2}\sum^k_{j=1}\alpha_{ij}(\zeta){\sum^n_{\nu=1}{\partial{\tilde\rho_j}\over\partial\overline\zeta_{\nu}}{\partial\over\partial\zeta_{\nu}}}-{1\over 2}\sum^k_{j=1}\beta_{ij}(\zeta){\partial\over\partial\overline\zeta_{\nu_j}}$$
where $\Big[\alpha_{ij}(\zeta)\Big]=A^{-1}$ and $\Big[\beta_{ij}(\zeta)\Big]=B^{-1}$.

\noindent Now it is easy to check that

$~~Y^{\zeta}_i\phi_j(\zeta,\zeta)=\delta_{ij}\mbox {~and~}Y^{\zeta}_i\tilde\rho_j=0\mbox {~for ~all ~}1\le i,j\le k.\hfill\Box$

\vspace*{0.5cm}

Let us introduce the following class of kernels for $\delta\ge 0$,
$${\cal L}_\delta={{\cal E}^j\over \prod^{k+1}_{i=1}(\phi_i+\delta)^{r_i}\
}$$
where $$2n-1-2\sum^{k+1}_{i=1}{r_i}+j\ge 0$$
and $$r_i\ge 1 \hbox{~for~all~}1\le i\le k+1$$
\paragraph{{\sc Remark 3.2.}}{\em Notice that the kernel $R$ is a finite sum of kernels  of type ${\cal L}_0$, and estimate (2.23) with estimate (2.24) hold, independently of $\delta$, for kernels ${\cal L}_\delta$.}

\vspace*{0.5cm}

If we denote by $X^z$ a tangential vector field to $M$ in $z$-variable and $X^\zeta$ the corresponding operator in $\zeta$-coordinates , then we have the following 
\paragraph{{\sc Lemma 3.3}}{\em Let $\delta>0$, then we have
$$X^z{\cal L}_{\delta}=-X^{\zeta}{\cal L}_{\delta}+\sum^k_{i=1}{(X^z+X^{\zeta})\phi_i\over Y^{\zeta}_i\phi_i}Y^{\zeta}_i({\cal L}_{\delta})+S_{\delta}.$$
where $S_{\delta}$ is a finite sum of kernels of type ${\cal L}_{\delta}$.}

\vspace*{0.5cm}

{\sl Proof.} It is not difficult to see that the following facts are true:

${\displaystyle {~(i)~}(X^z+X^{\zeta}){\cal E}^j \mbox {~ is~ of~ type ~}{\cal E}^j.}$

${\displaystyle{~(ii)~}(X^z+X^{\zeta})\phi_i \mbox {~is ~of~ type~}{\cal E}^1.}$

${\displaystyle{~(iii)~}|Y^{\zeta}_i\phi_i(\zeta,z)|\ge C \mbox
  {~for~} |\zeta-z|\le \epsilon\mbox, \varepsilon > 0,  {~and~} 1\le i\le k \mbox {~(see Lemma ~3.1)~}}$

${\displaystyle~(i\nu)\mbox{~If~} i\ne j \mbox {~then~} Y^{\zeta}_i\phi_j \mbox {~is~of~type~} {\cal E}^1 \mbox {~(cf.~Lemma ~3.1)~}.}$

${\displaystyle{~(\nu )~} Y^{\zeta}_i\phi_{k+1}-y_iY^{\zeta}_i\phi_i \mbox\
 {~is~of~type~} {\cal E}^1{~for~}1\le i\le k\mbox {~(see (2.25)}, (i\nu)\mbox{, Lemma ~3.1)}}$

${\displaystyle{~(\nu
    i)~}(X^z+X^{\zeta})\phi_{k+1}-{\sum^{k}_{i=1}y_i(X^z+X^{\zeta}) \phi_i}\mbox{~is~of~type~}{\cal E}^2 \mbox {~(see (2.25)~and~(i)~).}}$

${\displaystyle{~(\nu ii)~}(X^z+X^{\zeta})({1\over {\phi^{r_i}_i}})={(X^z+X^{\zeta})\phi_i\over Y^{\zeta}_i\phi_i}Y^{\zeta}_i({1\over {\phi^{r_i}_i}}).}$

\noindent The lemma follows now by a straightforward computation.$\hfill\Box$

\vspace*{0.5cm}

Now let $ \Omega\sub M_0$ and $f\in L^{\infty}(\Omega)\cap{{\cal C}_\star}^\ell (\Omega)$. Let $z_1\in\Omega$ and $\chi$ a smooth compactly supported  function on $\Omega$ such that$$\chi(\zeta)   =\cases{
0 & if $|\zeta-z_1|\ge{\epsilon\over 2}$\cr
1 & if $|\zeta-z_1|\le{\epsilon\over 4}$.\cr}$$
where $\epsilon $ is chosen so that (see Lemma 3.1) 
$$|Y^{\zeta}_i\phi_i(\zeta,z)|\ge C \mbox {~~for~} |\zeta-z|\le \epsilon\mbox {~and~ all~} 1\le i\le k.$$
Set $K:=\{z\in\Omega /|z-z_1|\le {\epsilon\over 4}\}$.

\noindent We write 
$$\int_\Omega f(\zeta)\wedge R(\zeta,z)=\int_\Omega \chi(\zeta) f(\zeta)\wedge R(\zeta,z)+\int_\Omega (1-\chi(\zeta)) f(\zeta)\wedge R(\zeta,z).$$
Let $J_1(f)$ denote the first integral in the right-hand side and $J_2(f)$ the second one.

\noindent Since $R(\zeta,z)$ is of class ${\cal C}^\infty$ in $z$ for $\zeta\ne z$ then $J_2(f)$ is of class ${\cal C}^\infty $ on $K$.

\noindent By Remark 3.2 to estimate $J_1(f)$, it is enough to do so for $\int_\Omega \chi f\wedge {\cal L}_0(\cdot,z)$.

\vspace*{0.5cm}

We have 
$$\int_\Omega\chi f\wedge {\cal L}_0(\cdot,z)=\lim_{\delta\to 0}\int_\Omega\chi f\wedge {\cal L}_{\delta}(\cdot,z).$$
By Lemma 3.3, we obtain from Stokes' theorem
\begin{eqnarray*}
X^z\int_\Omega\chi f\wedge {\cal L}_{\delta}(\cdot,z)&=&\pm \int_\Omega X^{\zeta}(\chi f)\wedge {\cal L}_{\delta}(\cdot,z)\\
&\pm& \sum^k_{i=1} \int_\Omega Y^{\zeta}_i(\chi f)\wedge {(X^z+X^{\zeta})\phi_i\over Y^{\zeta}_i\phi_i}{\cal L}_{\delta}(\cdot,z) + \int_\Omega \chi f \wedge S_{\delta}(\cdot,z).
\end{eqnarray*}
where $S_{\delta}$ is a finite sum of kernels of type ${\cal L}_{\delta}$.

\noindent Now , if we apply $r\le\ell$ derivatives , we can write 
$$X^z_1\cdots X^z_r \int_\Omega\chi f\wedge {\cal L}_{\delta}(\cdot,z)$$
as a sum of terms
$$\int_\Omega \tilde X^{\zeta}_1\cdots\tilde X^{\zeta}_j(\chi f)\wedge {\cal L}_{\delta}(\cdot,z)$$
with $0\le j\le r$.

\noindent Since 
$$\|\int_\Omega \tilde X^{\zeta}_1\cdots\tilde X^{\zeta}_j(\chi f)\wedge {\cal L}_{\delta}(\cdot,z)\|_{1\over 2}\le C \|f\|_{{\cal C}^\ell}$$
for $0\le j\le \ell$, independently of $\delta$ (see Remark 3.2), we conclude that
$$\int_\Omega \chi f\wedge {\cal L}_0(\cdot,z)\mbox{~ is~ of~ class~} {\cal C}^{\ell+{1\over 2}}\mbox{~ on~} \Omega.$$
Thus $J_1(f)$ is of class ${\cal C}^{\ell+{1\over 2}}$ on $\Omega$, and therefore 
$$ \int_\Omega f(\zeta)\wedge R(\zeta,z) \mbox {~is~ of~ class~} {\cal C}^{\ell+{1\over 2}}\mbox {~on~} K.$$
By noticing that $Y_i^z\Phi_j=-Y_i^{\zeta}\Phi_j+{\cal E}^1{~for~}1\le i,j\le k$ one  can show  in the same way 
$$ \int_\Omega f(z)\wedge R(\zeta,z) \mbox {~is~ of~ class~} {\cal C}^{\ell+{1\over 2}}\mbox {~on~} K$$ provided M is of class ${\cal C}^{\ell+3}$ (see (2.24)).
This completes the proof of the second part of Theorem 0.1 (cf.
Definition 2.9).$\hfill\Box$

\address{Institut Fourier\\
UFR de Math\'ematiques\\
UMR 5582 CNRS\\
 B.P 74 38402 Saint Martin d'H\`eres France\\
{\bf e-mail}: ybarkat@ujf-grenoble.fr}

\end{document}